\documentclass{article}

\usepackage[T2A]{fontenc}
\usepackage[cp1251]{inputenc}
\usepackage[english,russian]{babel}
\usepackage[tbtags]{amsmath}
\usepackage{amsfonts,amssymb}

\usepackage{mathrsfs}

\usepackage[paper]{mz2ewin}

\usepackage{graphicx}
\usepackage{url}
\numberwithin{equation}{section}

\theoremstyle{plain}

\newtheorem{theorem}{Теорема}

\newtheorem{lemma}{Лемма}

\newtheorem{definition}{Определение}

\theoremstyle{definition}

\newtheorem{remark}{Замечание}

\def\R{\mathbb R}
\def\E{\mathbb E}
\def\PP{\mathbb P}

\begin{document}

\udk{519.85}

\date{}

\author{Е.\,Л.~Гладин}
\address{Humboldt University of Berlin; \\
	Московский физико-технический институт, Москва; \\
	Институт проблем передачи информации РАН}
\email{egor.gladin@student.hu-berlin.de}

\author{А.\,В.~Гасников}
\address{Московский физико-технический институт, Москва; \\
	Институт проблем передачи информации РАН; \\ Кавказский Математический Центр Адыгейского государственного университета}
\email{gasnikov@yandex.ru}

\author{Е.\,С.~Ермакова}
\address{Московский физико-технический институт, Москва}
\email{yelena.ermakova@phystech.edu}

\title{Метод Вайды для задач выпуклой стохастической
	оптимизации небольшой размерности}

\markboth{Е.\,Л.~Гладин, А.\,В.~Гасников, Е.\,С.~Ермакова}{Метод Вайды для стохастической
	оптимизации}

\maketitle

\begin{fulltext}
\begin{abstract}
В работе рассматривается общая задача выпуклой стохастической оптимизации в пространстве небольшой размерности (например, 100 переменных). Известно, что для детерминированных задач выпуклой оптимизации небольших размеров наилучшим образом сходятся методы типа центров тяжести (например, метод Вайды). Для задач стохастической оптимизации вопрос о возможности использования метода Вайды сводится к вопросу о том, как он накапливает неточность в субградиенте. Недавний результат авторов об отсутствии накопления неточности на итерациях метода Вайды позволяет предложить его аналог для задач стохастической оптимизации. Основным приемом является замена субградиента в методе Вайды его пробатченным аналогом (средним арифметическим стохастических субградиентов). В настоящей работе осуществляется описанный план, что приводит в итоге к эффективному (в условиях возможности производить вычисления параллельно при батчинге) методу решения задач выпуклой стохастической оптимизации в пространствах небольших размеров.

Библиография: 15 названий.

\end{abstract}


\footnotetext{
	Работа Е.\,Л.~Гладина финансирована Deutsche Forschungsgemeinschaft (DFG, Немецкий исследовательский фонд) в рамках Excellence Strategy Германии – Берлинский математический исследовательский центр MATH+ (EXC-2046/1, project ID: 390685689).
	Работа А.\,В. Гасникова была выполнена при поддержке Министерства науки и высшего образования Российской Федерации (госзадание) № 075-00337-20-03, номер проекта 0714-2020-0005.
}

\section{Введение} \label{section_intro}
В работе рассматривается задача выпуклой стохастической оптимизации на компактном множестве простой структуры (что понимается под простым множеством поясняется далее в замечании \ref{simple_set}) в пространстве небольшой размерности $n$. Требуется решить задачу с точностью $\varepsilon$ по функции. Гладкость целевого функционала не предполагается. Известно, что в условиях отсутствия гладкости, возможность параллельно вычислять стохастические субградиенты целевого функционала (при разных реализациях случайности и в разных точках) слабо влияет на общее число последовательных итераций, которые требуется осуществить, чтобы достичь желаемой точности \cite{srebro2018}, \cite{bubeck2019}, \cite{diakonikolas2020}. Более точно, в условиях, при которых на каждой итерации разрешается вычислять не более $O\left(\text{Poly}(n,\varepsilon^{-1})\right)$ стохастических субградиентов, число последовательных итераций оптимальных методов будет (приводим зависимость только от $n$ и $\varepsilon$, опуская все остальные константы и логарифмические множители)
$$\min\left\{n,\frac{n^{1/3}}{\varepsilon^{2/3}},\frac{n^{1/4}}{\varepsilon},\frac{1}{\varepsilon^2}\right\}.$$
Более того, в предположении, что все запрошенные на любой итерации стохастические субградиенты могут быть посчитаны только в одной точке (определенной на текущей итерации) при этом (также как и раньше) с разными (стохастически независимыми) реализациями, число последовательных итераций оптимальных методов будет 
$$\min\left\{n,\frac{1}{\varepsilon^2}\right\}.$$
В данной работе нас будет интересовать последняя оценка и условия, при которых она получена.
Собственно, второй аргумент минимума отвечает ситуации, когда размерность пространства $n$ достаточно большая в сравнении с $\varepsilon^{-2}$. В этом случае параллелизация невозможна. Оптимальным методом будет обычный стохастический градиентный спуск.
В данной работе нас будет интересовать ситуация, когда размерность пространства небольшая. Тогда оптимальное число итераций будет $\widetilde{O}(n)$ ($\widetilde{O}(\cdot)$ означает то же самое, что и $O(\cdot)$ с точностью до логарифмических множителей по $\varepsilon$ и $n$). 
В этом случае параллелизация неизбежна, потому что нижняя оценка на число вызовов стохастического субградиента, чтобы достичь точности решения $\varepsilon$, даже для одномерных задач ($n = 1$) будет $\varepsilon^{-2}$ \cite{nemirovski1979}. Таким образом, даже при самом оптимистичном сценарии, число вычислений  стохастических субградиентов на каждой итерации будет не меньше, чем $\widetilde{O}\left(n^{-1}\varepsilon^{-2}\right)$. 
Однако, не известно такого метода, который бы так работал. Наилучший, известный нам на данный момент, метод из работы \cite{gladin2020ellipsoid} требует $\widetilde{O}\left(\varepsilon^{-2}\right)$ (параллельных) вычислений стохастических субградиентов на каждой итерации, и при этом нуждается в $\widetilde{O}(n^2)$ последовательных итерациях.
В настоящей статье будет описан метод, который также требует $\widetilde{O}\left(\varepsilon^{-2}\right)$ (параллельных) вычислений стохастических  субградиентов на каждой итерации (далее в статье мы будем называть эту величину размером батча), и при этом нуждается лишь в $\widetilde{O}(n)$ последовательных итерациях, что отвечает приведенной выше нижней оценке на число последовательных итераций.
Вопрос о возможности улучшения оценки на размер батча $\widetilde{O}\left(\varepsilon^{-2}\right)$ вплоть до $\widetilde{O}\left(n^{-1}\varepsilon^{-2}\right)$ остается открытым. По-видимому, такое улучшение невозможно.

В основу предлагаемого подхода положена конструкция минибатчинга \cite{dekel2012}, \cite{dvinskikh2020} и оптимальный по числу итераций (вычислений субградиента \cite{nemirovski1979}) метод (Вайды) решения малоразмерных задач выпуклой оптимизации из работы \cite{vaidya1989}, \cite{Bubeck}. Отметим, что метод Вайды имеет немного более дорогую теоретическую стоимость итерации (при достаточно больших $n$), чем более современный метод Ли--Сидфорда \cite{lee2014}, также оптимальный по числу итераций. Однако последний довольно сложно эффективно реализовать на практике. К тому же, в теории он становится лучше лишь при больших $n$, что не отвечает предположениям данной статьи.

\section{Описание результата} \label{section_1}
Рассмотрим задачу
\begin{equation}\label{problem:min_f}
	\min_{x \in Q} \left\{ f(x) := \E_{\xi} f(x, \xi) \right\},
\end{equation}
где $Q \subseteq \R^n$~--- выпуклое компактное множество с непустой внутренностью, размерность $n$ относительно небольшая (до ста), $f(x)$~--- непрерывная выпуклая функция. Будем считать, что стохастический субградиент $\partial_x f (x, \xi)$ удовлетворяет для некоторого $\sigma > 0$ условию (здесь и далее норма $\|\cdot\|$ считается евклидовой)
\begin{equation}
	\mathbb{E_\xi} \exp \left( \frac{\left\|\partial_x f (x, \xi) - g \right\|^2}{\sigma^2} \right) \leq \exp (1),
\end{equation}
где $g$~--- некоторый субградиент $f$ в точке $x$ (обозначение: $g \in \partial f (x)$).
\begin{definition}
	Случайный вектор $\tilde{x}$ со значениями в $Q$ будем называть $(\varepsilon, \beta)$-решением задачи \eqref{problem:min_f} для $\varepsilon>0,\ \beta \in (0,1)$, если $\PP \left\{  f(\tilde{x}) - f_* > \varepsilon \right\} \leq \beta$, где $f_*$~--- оптимальное значение в задаче \eqref{problem:min_f}.
\end{definition}
Такое определение позволяет ограничить вероятность того, что значение функции в приближённом решении сильно отклоняется от оптимального.

Пусть $\rho$ и $R$~--- радиусы некоторых евклидовых шаров $B_{\rho}$ и $B_R$ таких, что $B_{\rho} \subseteq Q \subseteq B_R$. Обозначим $B := \sup_{x, y \in Q} \lvert f(x) - f(y) \rvert$.
% Здесь и далее $\| \cdot \|$ означает евклидову норму.
В настоящей статье предлагается использовать метод Вайды с минибатчингом и доказывается, что он позволяет найти $(\varepsilon, \beta)$-решение задачи \eqref{problem:min_f} для $\varepsilon>0, \beta \in (0,1)$ за
\begin{equation*}
	N = O \left( n \log \left( \frac{n B R}{\rho \varepsilon} \right) \right)
\end{equation*}
итераций при размере батча\footnote{с точностью до множителя $\log \log \frac{n BR}{\rho \varepsilon}$}
\begin{equation*}
	r = O \left( \frac{\sigma^2 R^2}{\varepsilon^2} \log \frac{n}{\beta} \right).
\end{equation*}

\section{Метод Вайды} \label{section_2}
Метод Вайды (метод секущей плоскости) был предложен Вайдой в \cite{vaidya1989} для решения условной задачи оптимизации вида
\begin{equation}\label{problem_vaidya}
	\min_{x \in Q} f(x) ,
\end{equation}
где $Q \subset \R^n$~--- выпуклое компактное множество с непустой внутренностью, а целевая функция $f: Q \to \R$, определённая на $Q$, непрерывна и выпукла.

Далее вводятся обозначения и описывается алгоритм. Пусть $P = \{x \in \R^n: \, Ax\geq b\}$~--- ограниченный $n$-мерный многогранник, где  $A \in \R^{m\times n}$ и $b \in \R^m$. Логарифмический барьер множества $P$ определяется как
$$
L(x) = -\sum_{i=1}^{m} \log \left(a_{i}^{\top} x-b_{i}\right),
$$
где  $a_{i}^{\top}$~--- $i$-я строка матрицы $A$. Гессиан  $H(x)$  функции $L(x)$ равен
$$
H(x) =\sum_{i=1}^{m} \frac{a_{i} a_{i}^{\top}}{\left(a_{i}^{\top} x-b_{i}\right)^{2}}
$$
Матрица $H(x)$ положительно определена для всех $x$ из внутренности $P$. Волюметрический барьер (volumetric barrier) множества $P$ определяется как
$$
F(x) = \frac{1}{2} \log \left(\operatorname{det}H(x)\right),
$$
где $\operatorname{det}H(x)$ обозначает детерминант $H(x)$. Будем называть точку минимума функции $F$ на $P$ волюметрическим центром множества $P$. Обозначим за $\sigma_{i}(x)$ величины
\begin{equation}\label{volumetric_barrier}
	\sigma_{i}(x)=\frac{a_{i}^{\top} \left(H(x)\right)^{-1} a_{i}}{\left(a_{i}^{\top} x-b_{i}\right)^{2}}, \quad 1 \leq i \leq m,
\end{equation}

Параметры метода Вайды $\eta>0$ и $\gamma>0$~--- небольшие константы, такие что $\eta \leqslant 10^{-4}$ и $\gamma \leqslant 10^{-3} \eta$. Алгоритм производит последовательность пар $\left(A_k, b_k\right) \in \R^{m_k\times n}\times \R^{m_k}$ таких, что соответствующие многогранники содержат решение. В качестве начального многогранника, задаваемого парой $\left(A_0, b_0\right)$, обычно берётся симплекс (алгоритм может начинать с любого выпуклого ограниченного $n$-мерного многогранника, для которого легко вычислить волюметрический центр).
% ~--- например, с $n$-прямоугольника

В начале каждой итерации $k\geq 0$ находится приближённый волюметрический центр $x_k$ и вычисляются величины $\left\{\sigma_i(x_k)\right\}_{1\leq i \leq m}$. Следующий многогранник, характеризуемый парой $\left(A_{k+1}, b_{k+1}\right)$, получается из текущего в результате либо присоединения, либо удаления ограничения:
\begin{enumerate}
	\item Если для некоторого $i \in \{1,\ldots, m_k\}$ выполняется $\sigma_i(x_k) = \min\limits_{1\leq j\leq m}\sigma_j(x_k) < \gamma$, тогда $\left(A_{k+1}, b_{k+1}\right)$ получается исключением $i$-й строки из $\left(A_k, b_k\right)$.
	\item Иначе (если $\min\limits_{1\leq j\leq m}\sigma_j(x_k) \geq \gamma$) возможны два случая: 
	\begin{itemize}
		\item Если $x_k \notin Q$, то оракул возвращает вектор $c_k$ такой, что
		\begin{equation}
			\label{separation_oracle}
			c_k^{\top} (x - x_k) \geq 0\; \forall x \in Q;
		\end{equation}
		\item Если $x_k \in Q$, то оракул возвращает вектор $c_k$ такой, что
		$$f(x) \leq f(x_k)\; \forall x \in \left\{z \in Q: c_k^{\top} z\geq c_k^{\top} x_k \right\}, \text{ т.е. } c_k \in -\partial f(x_k).$$
	\end{itemize}
	Выберем $\beta_k \in \R$ таким, что
	$$
	\frac{c_k^{\top} \left(H(x_k)\right)^{-1} c_k}{\left(c_k^{\top} x_{k} -\beta_{k}\right)^{2}}=\frac{1}{2} \sqrt{\eta \gamma}.
	$$
	Определим $\left(A_{k+1}, b_{k+1}\right)$ добавлением строки $\left(c_k, \beta_k\right)$ к $\left(A_k, b_k\right)$.
\end{enumerate}
После $N$ итераций метод возвращает точку $\displaystyle x^N := \arg \min_{1 \leq k \leq N} f(x_k)$

Волюметрический барьер $F_k$ является самосогласованной функцией, поэтому может быть эффективно минимизирован методом Ньютона. Подробности и анализ метода Вайды можно найти в статье \cite{vaidya1989} и книге \cite{Bubeck}.
\begin{remark}\label{simple_set}
	Предполагается, что множейство $Q$ является простым в смысле, что можно проверить, что данная точка принадлежит $Q$, или вычислить вектор $c_k$, удовлетворяющий \eqref{separation_oracle}, за время, по порядку не превышающее время приближённого вычисления волюметрического центра. Поскольку для нахождения волюметрического центра достаточно нескольких шагов метода Ньютона, сложность такой операции характеризуется стоимостью обращения матрицы размера $n \times n$.
\end{remark}

Введём понятие неточного субградиента.
\begin{definition}\label{delta_subgradient}
	Пусть $\delta \geq 0,\ Q \subseteq \R^n$~--- выпуклое множество, $f: Q \to \R$~--- выпуклая функция. Вектор $g \in \R^n$ называется $\delta$-субградиентом $f$ в точке $x \in Q$, если
	\begin{equation*}
		f(y) \geq f(x) + \langle g, y-x \rangle - \delta\quad \forall y \in Q.
	\end{equation*}
	Множество $\delta$-субградиентов $f$ в точке $x$ обозначается $\partial_\delta f(x)$.
\end{definition}
Доказано \cite{gladin2021solving}, что в методе Вайды можно использовать $\delta$-субградиент вместо точного субградиента. А именно, справедлива следующая теорема:
\begin{theorem}{\cite{gladin2021solving}}\label{th:vaidya}
	Пусть $\rho$ и $R$~--- радиусы некоторых евклидовых шаров $B_{\rho}$ и $B_R$ таких, что $B_{\rho} \subseteq Q \subseteq B_R$, и пусть число $B>0$ таково, что $\left| f(x)-f(y) \right| \leq B\ \forall x,y \in Q$.
	После $N \geq \frac{2n}{\gamma} \log \left( \frac{n^{1.5} R}{\gamma \rho} \right) + \frac{1}{\gamma} \log \pi$ итераций метод Вайды с $\delta$-субгради"-ентом для задачи \eqref{problem_vaidya} возвращает такую точку $x^N$, что
	\begin{equation}
		f(x^N) - f(x_*) \leqslant \frac{n^{1.5} B R}{\gamma \rho} \exp \left( \frac{\log \pi -\gamma N}{2n} \right) + \delta,
	\end{equation}
	где $\gamma>0$~--- параметр алгоритма, $x_*$~--- решение задачи \eqref{problem_vaidya}.
\end{theorem}

\begin{remark}\label{mat_inversion}
	Помимо вычисления субградиента, в стоимость итерации метода Вайды входит стоимость обращения матрицы размера $n \times n$, что накладывает ограничение на его применение к задачам большой размерности.
\end{remark}

\section{Доказательство основного результата}
\label{section_3}
Введём следующее обозначение для стохастического субградиента с минибатчингом:
\begin{equation*} 
	\overset{r}{\partial_x} f \left( {x, \{ \xi^l \}_{l=1}^r} \right) := \frac{1}{r} \sum_{l=1}^r \partial_x f (x,  \xi^l ),
\end{equation*}
где $\xi^l$~---
% независимые реализации случайной величины $\xi$.
независимые случайные величины, распределённые одинаково с $\xi$.
Связать стохастический субградиент и $\delta$--субградиент позволяет следующий результат, полученный в работе \cite{gladin2020ellipsoid}.
\begin{lemma}{\cite{gladin2020ellipsoid}}\label{cons_batching}
	Пусть величина $\overset{r}{\partial_x} f \left( {x, \{ \xi^l \}_{l=1}^r} \right)$ была вычислена $N$ раз, $\tilde{g}_i$~--- её значение на шаге $i \in \overline{1, N}$ в точке $x^i$. Для любого $\beta \in (0, 1)$ справедливо
	\begin{equation*}
		\mathbb{P} \left( \bigcap\limits^N_{i=1} \left\{ \tilde{g}_i \in \partial_{\delta} f (x^i) \right\} \right) \geq 1 - \beta N, \quad \text{где } \delta = \left[\sqrt{2} + \sqrt{6\log{\beta^{-1}}}\right] \cdot \frac{\sigma R}{\sqrt{r}}.
	\end{equation*}
\end{lemma}
Таким образом, вероятность того, что величина $\overset{r}{\nabla}_x f \left( {x, \{ \xi^l \}_{l=1}^r} \right)$ будет являться $\delta$-субградиентом на всех $N$ шагах, составляет не менее $1 - \beta N$.

Основным результатом статьи является следующая теорема.
\begin{theorem}
	Метод Вайды с минибатчингом для задачи \eqref{problem:min_f} находит $(\varepsilon, \beta)$-ре"=шение для $\varepsilon>0, \beta \in (0,1)$ после
	\begin{equation*}
		N = \left\lceil \frac{2n}{\gamma} \log \left( \frac{n^{1.5} B R}{\gamma \rho \varepsilon} \right) + \frac{1}{\gamma} \log \pi \right\rceil
	\end{equation*}
	итераций при размере батча
	\begin{equation*}
		r = O \left( \frac{\sigma^2 R^2}{\varepsilon^2} \log \frac{n}{\beta} \right).
	\end{equation*}
	% АГ: ЧТО ТАКОЕ $D$? ОНО НЕ БЫЛО ОПРЕДЕЛЕНО.
\end{theorem}
\textit{Доказательство}\quad Согласно теореме \ref{th:vaidya},  $\varepsilon$-решение задачи \eqref{problem:min_f} может быть получено за
\begin{equation*}
	N = \left\lceil \frac{2n}{\gamma} \log \left( \frac{n^{1.5} B R}{\gamma \rho \varepsilon} \right) + \frac{1}{\gamma} \log \pi \right\rceil
\end{equation*}
итераций метода Вайды с $\frac{\varepsilon}{2}$-субградиентом. Воспользуемся леммой \ref{cons_batching} для определения размера батча $r$, необходимого для того, чтобы стохастический оракул возвращал $\frac{\varepsilon}{2}$-субградиент на каждой из $N$ итераций с вероятностью не менее $1 - \beta$:
\begin{equation*}
	\frac{\varepsilon}{2} = \left[\sqrt{2} + \sqrt{6\log \frac{N}{\beta}}\right] \cdot \frac{\sigma R}{\sqrt{r}} \Longrightarrow r = O \left( \frac{\sigma^2 R^2}{\varepsilon^2} \log \frac{n}{\beta} \right).
\end{equation*}

\section{Численный эксперимент}

Рассмотрим модель логистической регрессии для задачи классификации. Предсказанная вероятность принадлежности к классу с меткой 1 определяется по формуле
\begin{equation*}
	\hat{p}_x (w) = \frac{1}{1 + e^{-\langle w, x \rangle}},
\end{equation*}
где $x$~--- вектор признаков для объекта обучающей выборки (включая константный признак), $w$~--- веса модели. Обучающие примеры будем обозначать $\xi = (x, y)$, где $y \in \{0,1 \}$~--- метка класса объекта.
В качестве функции потерь выступает кросс-энтропия:
\begin{equation*}
	f_{\xi} (w) = y \ln \hat{p}_x (w) + (1-y) \ln \left(1 - \hat{p}_x (w) \right),
\end{equation*}
Задача оптимизации имеет вид
\begin{equation*}
	\min_{w \in Q} \left\{ f(w) := \E_{\xi} f_{\xi} (w) \right\},
\end{equation*}
где в качестве $Q$ можно взять евклидов шар достаточно большого радиуса. Таким образом, целевая функция является выпуклой и непрерывной, и минимизация осуществляется на компактном множестве с непустой внутренностью, что соответствует предположениям, в которых выведена оценка сложности предлагаемого метода.

В ходе эксперимента была использована выборка Covertype \cite{dataset}, состоящая из 581 012 объектов, 20\% из которых были выделены в тестовую выборку, а остальные 80\% использованы для обучения модели логистической регрессии с помощью метода Вайды и стохастического градиентного спуска. Количество признаков в данной выборке равно 55 (включая константный признак). Для сравнения методов измерено среднее значение ошибки на тестовой выборке на каждой итерации. 
Согласно выведенным оценкам скорости сходимости, метод Вайды требует большой размер батча $r$, однако в эксперименте он показывает хорошую сходимость уже при $r=128$. При увеличении размера батча вдвое сходимость становится ещё устойчивее, но большого ускорения не наблюдается. При $r=128$ SGD стабильно сходится при размере шага 0.1. При уменьшении батча вдвое кривая сходимости становится немного более зашумлённой, а при увеличении выше 128~-- никак не меняется, поэтому б\'oльшие размеры батча не отражены на графике. SGD показывает сходимость даже при небольших размерах батча, таких как $r=4$, но это требует уменьшения размера шага, что, в свою очередь, замедляет сходимость.
После каждой итерации обучения модели на отложенной тестовой выборке рассчитывалось значение функции потерь.
Зависимость среднего значения ошибки на тестовой выборке от количества итераций для каждого из методов отражена на рисунке \ref{paint1}. Кривая для метода Вайды при $r=256$ останавливается на моменте, где был совершён полный проход по обучающей выборке. Как видно из изображения, в данной задаче метод Вайды сходится значительно быстрее и приводит к меньшим значениям функции потерь.

\begin{figure}[t]
	\centering
	%     \begin{minipage}[h]{0.49\linewidth}
		%     \center{\includegraphics[width=0.9\linewidth]{ell_exp.eps}}
		% \end{minipage}
	% \includegraphics[scale=0.8,bb=00 00 425 270]{ell_exp.eps}
	\includegraphics[width=0.8\textwidth]{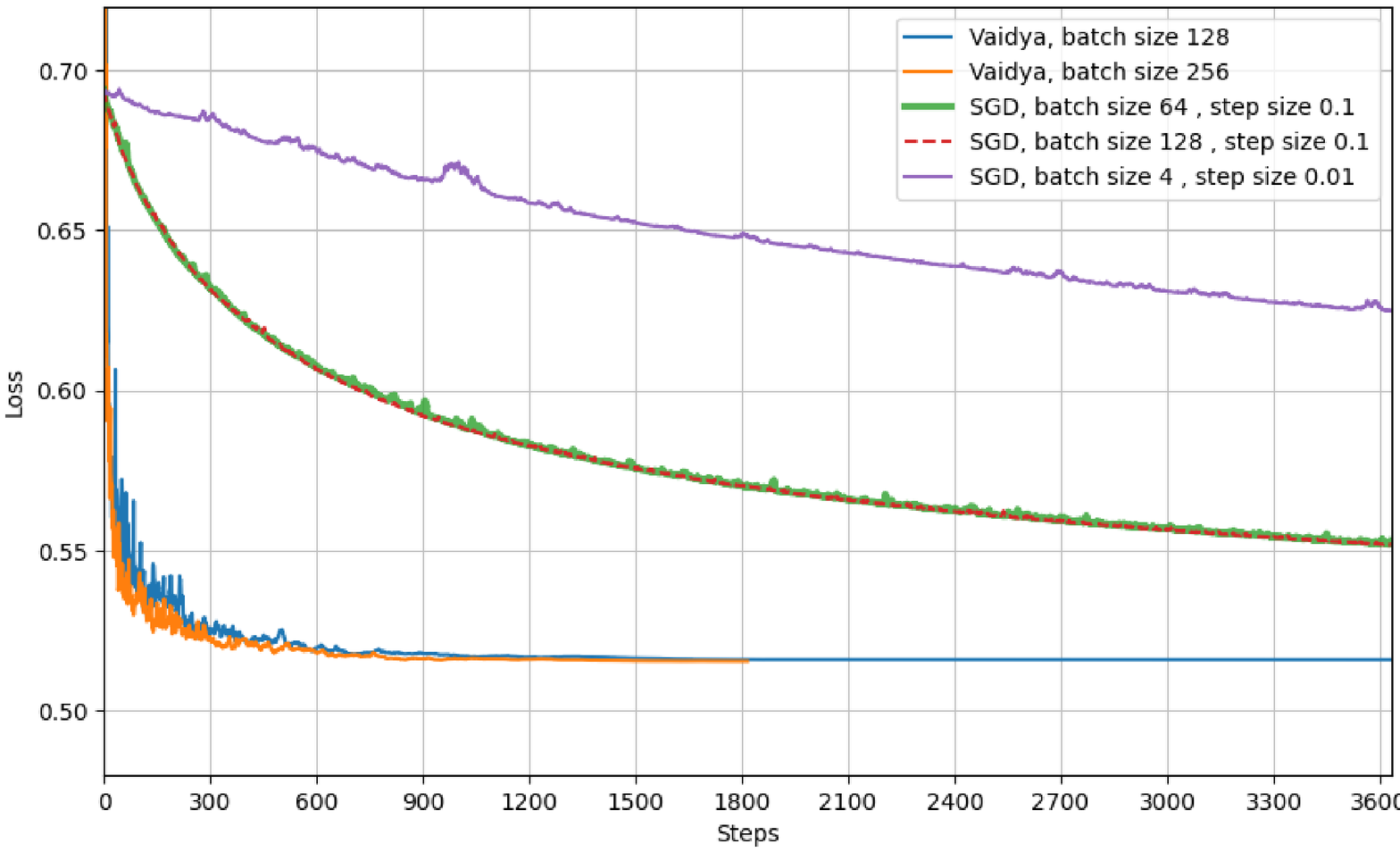}
	\caption{Сходимость метода Вайды и стохастического градиентного спуска при обучении логистической регрессии.}
	\label{paint1}
\end{figure}

Имплементацию метода Вайды и эксперименты можно найти в GitHub репозитории \cite{repo}.

\section{Заключение}
Используя результаты \cite{diakonikolas2020}, \cite{dekel2012},  приведенные в статье результаты могут быть распространены на случай, когда используется не евклидова норма. Вопрос выбора нормы сильно влияет на размер батча. Подходящий выбор может дополнительно улучшить зависимость размера батча от $n$ через константы $\sigma$ и $R$. 

Отметим также, что полученные в статье результаты представляют интерес не только в случае, когда рассматривается задача выпуклой негладкой стохастической  оптимизации, но и задача выпуклой (и даже сильно выпуклой) гладкой (функционал имеет липшицев градиент) стохастической оптимизации. Отмеченные выше возможности параллелизации вычислений при вычислении батча делают описанный подход конкурентным с различными вариантами ускоренных пробатченных методов (в том числе с редукцией дисперсии) и даже их тензорных аналогов. Детали см. в работе \cite{agafonov2021} и цитированной в ней литературе. 

Отметим в этой связи альтернативный (распределенный / оффлайн) способ получения основного результата данной статьи. Для этого заметим, что для получения решения задачи \eqref{problem:min_f} с точностью $\varepsilon$ нужно решить задачу 
$$\min_{x\in Q}\frac{1}{r}\sum_{k=1}^r f(x,\xi^k) + \frac{\varepsilon}{2R^2}\|x\|_2^2$$
c точностью $\sim\varepsilon^3$ \cite{shalev2009}.
Это можно сделать обычным (нестохастическим) вариантом метода Вайды за $\widetilde{O}(n)$ последовательных итераций, на каждой итерации параллельно считается субградиент всей функции -- $r$ вычислений субградиентов слагаемых. При этом общий объем требуемой выборки $\left\{\xi^k\right\}$ при таком подходе будет $r\sim\varepsilon^{-2}$, а не $\widetilde{O}(nr)\sim n\varepsilon^{-2}$, как в описанном в статье подходе. Однако, в статье описан онлайн подход, не требующий хранения в памяти всей выборки. В частности, в описанном в статье подходе не допускается вычисление субградиента функции $f(x,\xi^k)$ по $x$ при заданном $\xi^k$ более чем в одной точке $x$, в отличие от только что описанного подхода. Такое ограничение может быть связано, в том числе, с соображениями приватности.  
\end{fulltext}

%\begin{figure}[ht]
%\centering
%\includegraphics{1nek80-1}
%\caption{Бифуркационная диаграмма для сферы.}
%\label{fig1}
%\end{figure}

%=================Список литературы====================

\end{document}